\documentstyle[amssymb,amsmath,12pt]{article}

\newtheorem{th}{Theorem}[section]
\newtheorem{lem}[th]{Lemma}
\newtheorem{prop}[th]{Proposition}

\newtheorem{defn}[th]{Definition}
\newenvironment{defn-new}{\begin{defn} \em}{\end{defn}}
\newtheorem{rem}[th]{Remark}
\newenvironment{rem-new}{\begin{rem} \em}{\end{rem}}
\newtheorem{ex}[th]{Example}
\newenvironment{ex-new}{\begin{ex} \em}{\end{ex}}
\newtheorem{prob}[th]{Problem}
\newenvironment{prob-new}{\begin{prob} \em}{\end{prob}}

\newenvironment{notation-new}{\begin{rem} \em}{\end{rem}}

\newenvironment{agr-new}{\begin{rem} \em}{\end{rem}}

\makeatletter \@addtoreset{equation}{section} \makeatother
\begin{document}



\bigskip

\begin{center}
{\large {\bf Semi-parallelism of normal Jacobi operator}}

\medskip

{\large {\bf for Hopf hypersurfaces in complex two-plane Grassmannians}}

\bigskip \bigskip

Konstantina Panagiotidou and Mukut Mani Tripathi
\end{center}

\bigskip

{\bf Abstract.} It is proved the non-existence of Hopf hypersurfaces in $%
G_{2}({\Bbb C}^{m+2})$, $m\geq 3$, whose normal Jacobi operator is
semi-parallel, if the principal curvature of the Reeb vector field is
non-vanishing and the component of the Reeb vector field in the maximal
quaternionic subbundle ${\bf {\frak D}}$ or its orthogonal complement ${\bf {%
{\frak D^{\bot }}}}$ is invariant by the shape operator.

\medskip


\medskip


\section{Introduction\label{sect-intro}}

A complex two-plane Grassmannian $G_{2}(C^{m+2})$ is the set of all $2$%
-dimensional linear subspaces in $C^{m+2}$. It is a symmetric space and is
equipped with both a Kaehler structure $J$ and a quaternionic Kaehler
structure $J$ with a canonical local basis $\{J_{1},J_{2},J_{3}\}$, which
does not contain $J$. \medskip

Let $M$ be a real hypersurface in $G_{2}({\Bbb C}^{m+2})$, $N$ a unit normal
vector field of $M$ and $A$ the shape operator of $M$ with respect to $N$.
The Reeb vector field of $M$ is the structure vector field given by $\xi
=-JN $. Apart from the Reeb vector field, there are three more vector fields
given by $\xi _{\nu }=-J_{\nu }N$, $\nu =1,2,3$. Consequently, we have two
distributions on $M$ given by $[\xi ]={\rm Span}\{\xi \}$ and ${\frak D}%
^{\perp }={\rm Span}\{\xi _{1},\xi _{2},\xi _{3}\}$. We denote by ${\frak D}$
the orthogonal complement of the distribution ${\frak D}^{\perp }$ such that
$T_{p}M={\frak D}_{p}\oplus{\frak D}^{\perp }_{p}$, for each point $p \in M$.

\medskip

An important geometric condition for real hypersurfaces is the invariantness
of the distributions $[\xi ]$ and ${\frak D}^{\perp }$ under the action of
the shape operator. Under this condition, using a result due to Alekseevskii
\cite{Alekseevski-68}, Berndt and Suh classified the real hypersurfaces in
the following:

\begin{th}
\label{th-A} {\em (Theorem 1, \cite{Berndt-Suh-99})} Let $M$ be a connected
real hypersurface in $G_{2}({\Bbb C}^{m+2})$, $m\geq 3$. Then both the
distributions $\left[ \xi \right] $ and ${\frak D}^{\perp }$ are invariant
under the shape operator of $M$ if and only if either

\begin{itemize}
\item $M$ is of type {\bf (A)}, that is $M$ is an open part of a tube around
a totally geodesic $G_{2}\left( {\Bbb C}^{m+1}\right) $ in $G_{2}({\Bbb C}%
^{m+2})$, or

\item $M$ is of type {\bf (B)}, that is $m$ is even, say $m=2n$, and $M$ is
an open part of a tube around a totally geodesic ${\Bbb H}P^{n}$ in $%
G_{2}\left( {\Bbb C}^{2n+2}\right) $.
\end{itemize}
\end{th}

A real hypersurface $M$ in $G_{2}({\Bbb C}^{m+2})$ is said to be a Hopf
hypersurface if the Reeb vector field $\xi $ is principal, that is $A\xi
=\alpha \xi $, where $\alpha =g(A\xi ,\xi )$ is the corresponding principal
curvature to $\xi $. In such a case the integral curves of the Reeb vector
field $\xi $ are geodesics (Berndt and Suh \cite{Berndt-Suh-02}). Of course,
all of hypersurfaces in $G_{2}({\Bbb C}^{m+2})$ mentioned in Theorem~\ref%
{th-A} are Hopf hypersurfaces.

\medskip

In \cite{Berndt-91}, Berndt introduced the notion of {\em normal Jacobi
operator}
\[
\overline{R}_{N}(X)=\overline{R}(X,N)N\in {\rm End}(T_{x}M),\qquad x\in M,
\]%
for a real hypersurface $M$ in quaternionic projective spaces ${\Bbb H}P^{m}$
and in quaternionic hyperbolic spaces ${\Bbb H}H^{m}$, where $\overline{R}$
is the curvature tensor of the ambient space. He also proved the equivalence
of the commutation of $\overline{R}_{N}$ with the shape operator $A$ with
the fact that the distributions ${\frak D}$ and ${\frak D}^{\perp }$ are
invariant under the shape operator $A$.

\medskip

The classification of real hypersurfaces in $G_{2}({\Bbb C}^{m+2})$, whose
normal Jacobi operator $\overline{R}_{N}$ satisfies certain geometric
conditions, is one of great importance in the area of Differential Geometry.
In \cite{Perez-Jeong-Suh-2007}, Perez et. al. proved that ${\frak D^{\perp }}
$-invariant real hypersurfaces in $G_{2}({\Bbb C}^{m+2})$, whose normal
Jacobi operator commutes with both the structure tensor $\varphi $ and the
shape operator $A$ are locally congruent to one of type (A). Recently in
\cite{Jeong-Suh-Tripathi-2012}, Jeong, Suh and the second author considered
Hopf hypersurfaces in $G_{2}({\Bbb C}^{m+2})$ which satisfy the following
two commuting conditions
\[
\varphi A\overline{R}_{N}X=\overline{R}_{N}\varphi AX,\;\; X\in TM \qquad
{\rm and}\qquad A\varphi \varphi _{1}X=\varphi \varphi _{1}AX,\;\; X\in
{\frak D}^{\perp };
\]%
and proved that such real hypersurfaces are locally congruent to one of type
(A). The first condition is equivalent to $({\cal L}_{\xi}\overline{R}%
_{N})X=(\nabla_{\xi}\overline{R}_{N})X$.

\medskip

There are many interesting results concerning the non-existence of real
hypersurfaces in $G_{2}({\Bbb C}^{m+2})$ under certain geometric conditions
on the normal Jacobi operator. In \cite{Jeong-Suh-2008}, Jeong and Suh
examined cases of real hypersurfaces in $G_{2}({\Bbb C}^{m+2})$, when the
normal Jacobi operator is Lie $\xi $-parallel, that is ${\cal L}_{\xi}%
\overline{R}_{N}=0$. More precisely, they proved the non-existence of real
hypersurfaces in $G_{2}({\Bbb C}^{m+2})$ with ${\cal L}_{\xi}\overline{R}%
_{N}=0$ and one of the conditions $\xi \in {\frak D}^{\perp }$ and $\xi \in
{\frak D}$. They also proved the non-existence of Hopf hypersurfaces in $%
G_{2}({\Bbb C}^{m+2})$ with ${\cal L}_{\xi}\overline{R}_{N} = 0$ and
commuting shape operator on the distribution ${\frak D}^{\perp}$.

\medskip

In \cite{Jeong-Lee-Suh-2011}, it was proved that a Hopf hypersurface in $%
G_{2}({\Bbb C}^{m+2})$ does not exist if the normal Jacobi operator is Lie
parallel and the integral curves of ${\frak D}$- and ${\frak D}^{\perp }$-
components of the Reeb vector field are totally geodesic. In \cite%
{Machado-Perez-Suh-2011}, Machado et.~al.\ proved the non-existence of Hopf
hypersurfaces in $G_{2}({\Bbb C}^{m+2})$ whose normal Jacobi operator is of
Codazzi type (that is, $(\nabla_{X}\overline{R}_{N})Y=(\nabla_{Y}\overline{R}%
_{N})X$ for any $X,Y\in TM$) and ${\frak D}$- or ${\frak D}^{\bot}$%
-component of $\xi$ is invariant by the shape operator. In \cite%
{Jeong-Kim-Suh-2010}, Jeong et.~al.$\!$ proved the non-existence of Hopf
hypersurfaces in $G_{2}({\Bbb C}^{m+2})$ with parallel normal Jacobi
operator, that is $\nabla _{X}\overline{R}_{N}=0$. In \cite{Jeong-Suh-2011},
the non-existence of Hopf hypersurfaces in $G_{2}({\Bbb C}^{m+2})$ whose
normal Jacobi operator is $\left( [\xi ]\cup {\frak D}^{\perp }\right) $%
-parallel, which is a weaker condition then the previous one, was proved.

\medskip

A tensor field $P$ of type $(1,s)$ on a Riemannian manifold is said to be
{\em semi-parallel} if $R\cdot P=0$, where $R$ is the curvature tensor of
the manifold and acts as a derivation on $P$ \cite{CK}. In the geometry of
real hypersurfaces in complex space form the following results concerning
the semi-parallelism conditions have been proved. In \cite{PS1}, Perez and
Santos proved that there exist no real hypersurfaces in complex projective
space $CP^{n}$, $n\geq 3$, with semi-parallel structure Jacobi operator
(that is $R\cdot R_{\xi }=0$, where $R_{\xi} = R(\cdot ,\xi)\xi$ and $\xi$
is the structure vector field). Later, Cho and Kimura \cite{CK} generalized
this work and proved that there do not exist real hypersurfaces in complex
space forms equipped with semi-parallel structure Jacobi operator. Finally,
Niebergall and Ryan in \cite{NR1} studied real hypersurfaces in complex
space forms equipped with the semi-parallel shape operator $A$.

\medskip

Motivated by these studies the following question is raised naturally:

\begin{prob-new}
Do there exist real hypesurfaces in $G_{2}({\Bbb C}^{m+2})$, $m\geq 3$,
whose normal Jacobi operator, structure Jacobi operator or shape operator is
semi-parallel?
\end{prob-new}

In the present paper we give the answer partially and prove the following:

\begin{th}
\label{th-main} There does not exist any connected Hopf hypersurface $M$ in $%
G_{2}({\Bbb C}^{m+2})$, $m\geq 3$, equipped with semi-parallel normal Jacobi
operator, if $\alpha \neq 0$ and ${\frak D}$- or ${\frak D^{\bot }}$-
-component of the Reeb vector field $\xi $ is invariant by the shape
operator $A$.
\end{th}

The paper is organized as follows. In section~\ref{sect-Riem-Geom-of-CTPG},
we give a brief description of complex two plane Grassmanians. In section %
\ref{sect-real-hyp-in-CTPG} basic relations for real hypersurfaces in $G_{2}(%
{\Bbb C}^{m+2})$ are presented. Section~\ref{sect-key-lemmas} contains some
key results for further use. Finally, in section~\ref{sect-main-th}, we give
the proof of Theorem~\ref{th-main}.

\section{Riemannian Geometry of $G_{2}({\Bbb C}^{m+2})$ \label%
{sect-Riem-Geom-of-CTPG}}

The complex two-plane Grassmannian $G_{2}({\Bbb C}^{m+2})$ is the Grassmann
manifold of all complex $2$-dimensional linear subspaces in ${\Bbb C}^{m+2}$%
. The special unitary group $G=SU(m+2)$ acts transitively on $G_{2}({\Bbb C}%
^{m+2})$ with stabilizer isomorphic to $K=S(U(2)\times U(m))\subset G$. Thus
$G_{2}({\Bbb C}^{m+2})$ can be identified with the homogeneous space $G/K$,
which can be equipped with the unique analytic structure for which the
natural action of $G$ on $G_{2}({\Bbb C}^{m+2})$ becomes analytic. Denote by
${\frak g}$ and ${\frak l}$ the Lie algebra of $G$ and $K$, respectively.
Let ${\frak m}$ be the orthogonal complement of ${\frak l}$ in ${\frak g}$
with respect to the Cartan-Killing form $B$ of ${\frak g}$. Then ${\frak g}=%
{\frak l}\oplus {\frak m}$ is an $Ad(K)$-invariant reductive decomposition
of ${\frak g}$. we put $o=eK$ and identify $T_{o}G_{2}({\Bbb C}^{m+2})$ with
${\frak m}$ in the usual manner. Since $B$ is negative definite on ${\frak g}
$, therefore the restriction $(-B)|_{{\frak m}\times {\frak m}}$ yields a
positive definite inner product on ${\frak m}$. By $Ad(K)$-invariance of $B$
this inner product can be extended to a $G$-invariant Riemannian metric $g$
on $G_{2}({\Bbb C}^{m+2})$. In this manner $G_{2}({\Bbb C}^{m+2})$ becomes a
Riemannian homogeneous symmetric space. For computational reasons we
normalize the Riemannian metric $g$ such that the maximal sectional
curvature of $(G_{2}({\Bbb C}^{m+2}),g)$ becomes $8$.

\medskip

When $m=1$, $G_{2}({\Bbb C}^{3})$ is isometric to the $2$-dimensional
complex projective space ${\Bbb C}P^{2}$ with constant holomorphic sectional
curvature $8$. When $m=2$, the isomorphism ${\rm Spin}(6)\backsimeq SU(4)$
provides an isometry between $G_{2}({\Bbb C}^{4})$ and the real Grassmann
manifold $G_{2}^{+}({\Bbb R}^{6})$ of oriented $2$-dimensional linear
subspaces of ${\Bbb R}^{6}$. Therefore, we usually assume that $m\geq 3$.

\medskip

The Lie algebra ${\frak l}$ has the direct sum decomposition ${\frak l}=%
{\frak su}(m)\oplus {\frak su}(m)\oplus \Re $, where $\Re $ is the center of
${\frak l}$. Regarding ${\frak l}$ as the holonomy algebra of $G_{2}({\Bbb C}%
^{m+2})$, the center $\Re $ induces a Kaehler structure $J$ and the ${\frak %
su}(2)$-part induces a quaternionic Kaehler structure ${\frak J}$ on $G_{2}(%
{\Bbb C}^{m+2})$. If $J_{\nu }$ is any almost Hermitian structure in ${\frak %
J}$, then \thinspace $JJ_{\nu }=J_{\nu }J$, and $JJ_{\nu }$ is a symmetric
endomorphism with $(JJ_{\nu })^{2}=I$ and ${\rm tr}(JJ_{\nu })=0$.

\medskip

A canonical local basis $\{J_{1},J_{2},J_{3}\}$ of ${\frak J}$ consists of
three local almost Hermitian structures $J_{\nu }$ in ${\frak J}$ such that $%
J_{\nu }J_{\nu +1}=J_{\nu +2}=-J_{\nu +1}J_{\nu }$, where the index is taken
modulo $3$. Since ${\frak J}$ is parallel with respect to the Riemannian
connection $\overline{\nabla }$ of $(G_{2}({\Bbb C}^{m+2}),g)$, there exist
for any canonical local basis $J_{1}$, $J_{2}$, $J_{3}$ of ${\frak J}$ three
local \thinspace $1$-forms $q_{1}$, $q_{2}$, $q_{3}$, such that
\begin{equation}
\overline{\nabla }_{X}J_{\nu }=q_{\nu +2}(X)J_{\nu +1}-q_{\nu +1}(X)J_{\nu
+2}  \label{eq-G2-der}
\end{equation}%
for all vector fields $X$ on $G_{2}({\Bbb C}^{m+2})$.

\medskip

The Riemann curvature tensor $\overline{R}$ of $G_{2}({\Bbb C}^{m+2})$ is
locally given by \cite{Berndt-97}
\begin{eqnarray}
\overline{R}(X,Y)Z &=&g(Y,Z)X-g(X,Z)Y\frac{{}}{{}}  \nonumber \\
&&+\ g\left( JY,Z\right) JX-g\left( JX,Z\right) JY-2g\left( JX,Y\right) JZ
\nonumber \\
&&+\ \sum_{\nu =1}^{3}\left\{ g\left( J_{\nu }Y,Z\right) J_{\nu }X-g\left(
J_{\nu }X,Z\right) J_{\nu }Y-2g\left( J_{\nu }X,Y\right) J_{\nu }Z\right\}
\nonumber \\
&&+\ \sum_{\nu =1}^{3}\left\{ g\left( J_{\nu }JY,Z\right) J_{\nu }JX-g\left(
J_{\nu }JX,Z\right) J_{\nu }JY\right\}  \label{eq-G2-curvature}
\end{eqnarray}%
for all vector fields $X,Y,Z$ on $G_{2}({\Bbb C}^{m+2})$, where $\left\{
J_{1},J_{2},J_{3}\right\} $ is any canonical local basis of ${\frak J}$.
This expression involves the Riemannian curvature tensor of $S^{4m}$, ${\Bbb %
C}P^{2m}$ and ${\Bbb H}P^{m}$.

\section{Real hypersurfaces in $G_{2}({\Bbb C}^{m+2})$ \label%
{sect-real-hyp-in-CTPG}}

Let $M$ be a real hypersurface in $G_{2}({\Bbb C}^{m+2})$, that is a
hypersurface of $G_{2}(C^{m+2})$ with real codimension one. The induced
Riemannian metric on $M$ is denoted by $g$ and $\nabla $ denotes the induced
Riemannian connection of $\left( M,g\right) $. Let $N$ be a local unit
normal field of $M$ and $A$ the shape operator of $M$ with respect to $N$.

\medskip

Now let us put
\begin{equation}
JX=\varphi X+\eta (X)N,\qquad J_{\nu }X=\varphi _{\nu }X+\eta _{\nu }(X)N
\label{eq-JX}
\end{equation}%
for any tangent vector $X$ of a real hypersurface $M$ in $G_{2}({\Bbb C}%
^{m+2})$.

\medskip

The Kaehler structure $J$ of $G_{2}({\Bbb C}^{m+2})$ induces a local almost
contact metric structure $\left( \varphi ,\xi ,\eta ,g\right) $ on $M$ in
the following way
\[
\varphi ^{2}X=-X+\eta (X)\xi ,\;\eta (X)=1,\;\varphi \xi =0,\;\eta
(X)=g(x,\xi ).
\]

If $M$ is orientable then $\xi $ is globally defined and is the induced Reeb
vector field on $M$. Furthermore, let $\left\{ J_{1},J_{2},J_{3}\right\} $
be a canonical local basis of ${\frak J}$. Then each $J_{\nu }$ induces an
almost contact metric structure $\left( \varphi _{\nu },\xi _{\nu },\eta
_{\nu },g\right) $ on $M$. Locally, the orthogonal complement of the real
span of $\xi $ in $TM$ is denoted by ${\frak H}$ and the orthogonal
complement of the real span of $\xi _{1}$, $\xi _{2}$, $\xi _{3}$ in $TM$ is
denoted by ${\frak D}$.

\medskip

In view of (\ref{eq-G2-curvature}), the Gauss equation is given by
\begin{eqnarray}
R(X,Y)Z &=&g(Y,Z)X-g(X,Z)Y  \nonumber \\
&&+\ g\left( \varphi Y,Z\right) \varphi X-g\left( \varphi X,Z\right) \varphi
Y-2g\left( \varphi X,Y\right) \varphi Z  \nonumber \\
&&+\ \sum_{\nu =1}^{3}\left\{ g\left( \varphi _{\nu }Y,Z\right) \varphi
_{\nu }X-g\left( \varphi _{\nu }X,Z\right) \varphi _{\nu }Y-2g\left( \varphi
_{\nu }X,Y\right) \varphi _{\nu }Z\right\}  \nonumber \\
&&+\ \sum_{\nu =1}^{3}\left\{ g\left( \varphi _{\nu }\varphi Y,Z\right)
\varphi _{\nu }\varphi X-g\left( \varphi _{\nu }\varphi X,Z\right) \varphi
_{\nu }\varphi Y\right\}  \nonumber \\
&&-\ \sum_{\nu =1}^{3}\left\{ \eta (Y)\eta _{\nu }(Z)\varphi _{\nu }\varphi
X-\eta (X)\eta _{\nu }(Z)\varphi _{\nu }\varphi Y\right\}  \nonumber \\
&&-\ \sum_{\nu =1}^{3}\left\{ \eta (X)g\left( \varphi _{\nu }\varphi
Y,Z\right) -\eta (Y)g\left( \varphi _{\nu }\varphi X,Z\right) \right\} \xi
_{\nu }  \nonumber \\
&&+\ g\left( AY,Z\right) AX-g\left( AX,Z\right) AY  \label{eq-Gauss-real-hyp}
\end{eqnarray}%
where $R$ denotes the curvature tensor of the real hypersurface $M$ in $%
G_{2}({\Bbb C}^{m+2})$.

\medskip

It is straightforward to verify the following identities%
\begin{equation}
\begin{array}{c}
\varphi _{\nu }\xi _{\nu +1}=\xi _{\nu +2},\qquad \varphi _{\nu +1}\xi _{\nu
}=-\,\xi _{\nu +2},\medskip \\
\varphi \xi _{\nu }=\varphi _{\nu }\xi ,\qquad \eta _{\nu }\left( \varphi
X\right) =\eta \left( \varphi _{\nu }X\right) ,\medskip \\
\varphi _{\nu }\varphi _{\nu +1}X=\varphi _{\nu +2}X+\eta _{\nu +1}(X)\xi
_{\nu },\medskip \\
\varphi _{\nu +1}\varphi _{\nu }X=-\,\varphi _{\nu +2}X+\eta _{\nu }(X)\xi
_{\nu +1}.%
\end{array}
\label{eq-induced-hyp}
\end{equation}

\medskip

In view of (\ref{eq-JX}), (\ref{eq-G2-der}) and (\ref{eq-induced-hyp}), it
is known that
\[
\left( \nabla _{X}\varphi \right) Y=\eta (Y)AX-g\left( AX,Y\right) \xi
,\qquad \nabla _{X}\xi =\varphi AX,
\]%
\[
\nabla _{X}\xi _{\nu }=q_{\nu +2}(X)\xi _{\nu +1}-q_{\nu +1}(X)\xi _{\nu
+2}+\varphi _{\nu }AX,
\]%
\[
\left( \nabla _{X}\varphi _{\nu }\right) Y=-\,q_{\nu +1}(X)\varphi _{\nu
+2}Y+q_{\nu +2}(X)\varphi _{\nu +1}Y+\eta _{\nu }(Y)AX-g\left( AX,Y\right)
\xi _{\nu }.
\]%
Summing up these formulas, we also find the following
\begin{eqnarray*}
\nabla _{X}\left( \varphi _{\nu }\xi \right) &=&\left( \nabla _{X}\varphi
_{\nu }\right) \xi +\varphi _{\nu }\left( \nabla _{X}\xi \right) \\
&=&-\,q_{\nu +1}(X)\varphi _{\nu +2}\xi +q_{\nu +2}(X)\varphi _{\nu +1}\xi \\
&&+\,\eta _{\nu }\left( \xi \right) AX-g\left( AX,\xi \right) \xi _{\nu
}+\varphi _{\nu }\varphi AX.
\end{eqnarray*}%
Moreover, from $JJ_{\nu }=J_{\nu }J$, $\nu =1,2,3$, it follows that
\[
\varphi _{\nu }\varphi X=\varphi \varphi _{\nu }X-\eta _{\nu }(X)\xi +\eta
(X)\xi _{\nu }.
\]%
\noindent For more details we refer to \cite{Alekseevski-68}, \cite%
{Berndt-97}, \cite{Berndt-Suh-99} and \cite{Berndt-Suh-02}. \medskip

\section{Key Lemmas\label{sect-key-lemmas}}

We consider a connected, orientable, Hopf hypersurface $M$ in $%
G_{2}(C^{m+2}) $ with $\alpha \neq 0$ and semi-parallel normal Jacobi
operator. The normal Jacobi operator $\overline{R}_{N}$ for a real
hypersurface $M$ in $G_{2}(C^{m+2})$ is given by
\begin{eqnarray}
\overline{R}_{N}(X) &=&X+3\eta (X)\xi +3\sum_{\nu =1}^{3}\eta _{\nu }(X)\xi
_{\nu }  \nonumber \\
&&-\sum_{\nu =1}^{3}\{\eta _{\nu }(\xi )\left( \varphi _{\nu }\varphi X-\eta
(X)\xi _{\nu }\right) -\eta _{\nu }(\varphi X)\varphi _{\nu }\xi \}
\label{eq-normal-Jacobi-operator}
\end{eqnarray}%
for any vector field $X$ tangent to $M$. Furthermore, semi-parallelism
condition of it, that is $R(X,Y)\cdot \overline{R}_{N}=0$, implies
\begin{equation}
R(X,Y)\overline{R}_{N}Z=\overline{R}_{N}(R(X,Y)Z)  \label{eq-semi-parallel}
\end{equation}%
for all vector fields $X,Y,Z$ tangent to $M$.

\begin{lem}
\label{lem-1} Let $M$ be a Hopf hypersurface in $G_{2}({\Bbb C}^{m+2})$ such
that ${\frak D}$- or ${\frak D}^{\bot }$-component of $\xi $ is invariant by
the shape operator $A$ and $\alpha \neq 0$. If the normal Jacobi operator is
semi-parallel, then $\xi \in {\frak D}$ or $\xi \in {\frak D}^{\bot }$.
\end{lem}

\noindent {\bf Proof.} Suppose that $\xi $ is written as
\begin{equation}
\xi =\eta (U)U+\eta (\xi _{1})\xi _{1}+\eta (\xi _{2})\xi _{2}+\eta (\xi
_{3})\xi _{3},  \label{eq-xi-representation}
\end{equation}%
where $U$ is a unit vector in ${\frak D}$ and $\eta (U)\neq 0$ and $\eta
(\xi _{\kappa })\neq 0$ for at least one $\kappa \in \left\{ 1,2,3\right\} $%
. Then relation (\ref{eq-xi-representation}) implies that
\begin{equation}
\varphi _{\kappa }\xi =\eta (U)\varphi _{\kappa }U+\eta (\xi _{\kappa
+1})\xi _{\kappa +2}-\eta (\xi _{\kappa +2})\xi _{\kappa +1}.
\label{eq-phi-k-xi-representation}
\end{equation}%
From (\ref{eq-normal-Jacobi-operator}), we get
\begin{equation}
\overline{R}_{N}(\xi )=4\xi +4\sum_{\nu =1}^{3}\eta (\xi _{\nu })\xi _{\nu },
\label{eq-RN-xi}
\end{equation}%
\begin{equation}
\overline{R}_{N}(\xi _{\kappa })=4\xi _{\kappa }+4\eta (\xi _{\kappa })\xi
+2\eta (\xi _{\kappa +1})\varphi _{\kappa +2}\xi -2\eta (\xi _{\kappa
+2})\varphi _{\kappa +1}\xi ,  \label{eq-RN-xi-k}
\end{equation}%
\begin{equation}
\overline{R}_{N}(\varphi _{\kappa }\xi )=2\eta (\xi _{\kappa +1})\xi
_{\kappa +2}-2\eta (\xi _{\kappa +2})\xi _{\kappa +1}.
\label{eq-RN-phi-k-xi}
\end{equation}%
Since the normal Jacobi operator is semi-parallel, from (\ref%
{eq-semi-parallel}) and (\ref{eq-RN-xi}), we get
\begin{equation}
\overline{R}_{N}(R(\xi ,\xi _{\kappa })\xi ) = 4R(\xi ,\xi _{\kappa })\xi
+4\sum_{\nu =1}^{3}\eta (\xi _{\nu })R(\xi ,\xi _{\kappa })\xi _{\nu }.
\label{eq-RN(R(xi,xi-{kappa})xi)}
\end{equation}%
Since ${\frak D}$- or ${\frak D}^{\bot }$-component of $\xi $ is assumed to
be invariant by the shape operator $A$, we obtain
\begin{equation}
AU=\alpha U\quad {\rm and}\quad A\xi _{\kappa }=\alpha \xi _{\kappa },\qquad
\kappa \in \left\{ 1,2,3\right\} .  \label{eq-A-invariant}
\end{equation}%
In view of (\ref{eq-A-invariant}), from relation (\ref{eq-Gauss-real-hyp})
we get
\begin{equation}
R(\xi ,\xi _{\kappa })\xi =\alpha ^{2}\eta (\xi _{\kappa })\xi -\alpha
^{2}\xi _{\kappa }+2\eta (\xi _{\kappa +1})\varphi _{\kappa +2}\xi -2\eta
(\xi _{\kappa +2})\varphi _{\kappa +1}\xi .  \label{eq-R(xi,xi-{kappa})xi}
\end{equation}%
Substituting (\ref{eq-R(xi,xi-{kappa})xi}) in (\ref%
{eq-RN(R(xi,xi-{kappa})xi)}), we lead to the following
\begin{eqnarray}
4\sum_{\nu =1}^{3}\eta (\xi _{\nu })R(\xi ,\xi _{\kappa })\xi _{\nu }
&=&\alpha ^{2}\eta (\xi _{\kappa })\overline{R}_{N}(\xi )-\alpha ^{2}%
\overline{R}_{N}(\xi _{\kappa })  \nonumber \\
&&+2\eta (\xi _{\kappa +1})\overline{R}_{N}(\varphi _{\kappa +2}\xi )-2\eta
(\xi _{\kappa +2})\overline{R}_{N}(\varphi _{\kappa +1}\xi )  \label{eq-D7}
\\
&&-4\alpha ^{2}\eta (\xi _{\kappa })\xi +4\alpha ^{2}\xi _{\kappa }
\nonumber \\
&&-8\eta (\xi _{\kappa +1})\varphi _{\kappa +2}\xi +8\eta (\xi _{\kappa
+2})\varphi _{\kappa +1}\xi .  \nonumber
\end{eqnarray}%
Taking the inner product of (\ref{eq-D7}) with $U$, in view of (\ref%
{eq-RN-xi-k}), (\ref{eq-RN-phi-k-xi}) and (\ref{eq-phi-k-xi-representation})
we obtain
\begin{equation}
\sum_{\nu =1}^{3}\eta (\xi _{\nu })g(R(\xi ,\xi _{\kappa })\xi _{\nu
},U)=-\alpha ^{2}\eta (\xi _{\kappa })\eta (U).  \label{eq-D9}
\end{equation}%
We calculate $R(\xi ,\xi _{\kappa })\xi _{\nu }$ from relation (\ref%
{eq-Gauss-real-hyp}) taking into account (\ref{eq-A-invariant}) and then we
take the inner product with $U$ and we lead to the following relation
\begin{equation}
g(R(\xi ,\xi _{\kappa })\xi _{\nu },U)=\alpha ^{2}\eta _{\kappa }(\xi _{\nu
})\eta (U).  \label{eq-D10}
\end{equation}%
From (\ref{eq-D9}) and (\ref{eq-D10}) we get
\[
\alpha ^{2}\eta (\xi _{\kappa })\eta (U)=0,\qquad \kappa \in \left\{
1,2,3\right\} ,
\]%
which is a contradiction. $\blacksquare $

\medskip

Now, we examine the case when the Reeb vector field $\xi $ belongs to the
distribution ${\frak D}^{\perp }$. In fact, we have the following

\begin{lem}
\label{lem-2} Let $M$ be a Hopf hypersurface in $G_{2}({\Bbb C}^{m+2})$ and $%
\alpha \neq 0$, with semi-parallel normal Jacobi operator and $\xi \in
{\frak D}^{\bot }$ then $g(A{\frak D},{\frak D}^{\bot })=0$.
\end{lem}

\noindent {\bf Proof.} Let $W\in {\frak D}$ arbitrarily. In order to prove
that $g(A{\frak D},{\frak D}^{\bot })=0$, it suffices to prove that $%
g(AW,\xi _{\kappa })=0$, $\kappa =1,2,3$. Since $\xi \in {\frak D^{\bot }}$,
we have that $JN\in {\frak J}N$. Let $J_{1}$ be an almost Hermitian
structure of ${\frak J}$ such that $JN=J_{1}N$. Then we obtain that $\xi
=\xi _{1}$ and $\eta (\xi _{2})=\eta (\xi _{3})=0$. Furthermore, $\varphi
\xi _{2}=-\xi _{3}$, $\varphi \xi _{3}=\xi _{2}$ and $\varphi ({\frak D})
\subset {\frak D}$.

\medskip

Due to the fact that $M$ is a Hopf hypersurface, we have that $A\xi =\alpha
\xi $ and so $g(AW,\xi )=g(AW,\xi _{1})=0$. Thus, it remains to prove that
\[
g(AW,\xi _{\kappa })=0,\qquad \kappa =2,3{\frak .}
\]%
From (\ref{eq-normal-Jacobi-operator}), we obtain
\begin{equation}
\overline{R}_{N}(\xi )=8\xi ,\qquad \overline{R}_{N}(W)=W-\varphi
_{1}\varphi W.  \label{eq-f4}
\end{equation}%
Using (\ref{eq-f4}) in (\ref{eq-semi-parallel}) we get
\begin{equation}
8R(W,\xi )\xi =\overline{R}_{N}(R(W,\xi )\xi ).  \label{eq-f5}
\end{equation}%
In view of $A\xi =\alpha \xi $, from (\ref{eq-Gauss-real-hyp}), it follows
that
\begin{equation}
R(W,\xi )\xi =W+\alpha AW-\varphi _{1}\varphi W.  \label{eq-f6}
\end{equation}%
Substituting (\ref{eq-f6}) in (\ref{eq-f5}) and taking into consideration (%
\ref{eq-f4}) we lead to the following
\begin{equation}
8W+8\alpha AW-8\varphi _{1}\varphi W=\overline{R}_{N}(W)+\alpha \overline{R}%
_{N}(AW)-\overline{R}_{N}(\varphi _{1}\varphi W).  \label{eq-f7}
\end{equation}%
From (\ref{eq-normal-Jacobi-operator}) we also get
\[
\overline{R}_{N}(AW)=AW+2\eta _{2}(AW)\xi _{2}+2\eta _{3}(AW)\xi
_{3}-\varphi _{1}\varphi AW,
\]%
\[
\overline{R}_{N}(\varphi _{1}\varphi W)=\varphi _{1}\varphi W-\varphi
_{1}\varphi (\varphi _{1}\varphi W).
\]%
Substitution of the previous two relations in (\ref{eq-f7}) gives
\[
7W+7\alpha AW-6\varphi _{1}\varphi W=2\alpha \eta _{2}(AW)\xi _{2}+2\alpha
\eta _{3}(AW)\xi _{3}+\varphi _{1}\varphi (\varphi _{1}\varphi W)-\alpha
\varphi _{1}\varphi AW.
\]%
Taking the inner product of the last relation with $\xi _{\kappa }$, $\kappa
=2,3$, and because of $\alpha \neq 0$ implies
\[
\eta _{\kappa }(AW)=0,\qquad \kappa =2,3,
\]%
and this completes the proof. $\blacksquare $

\medskip

Finally, in the case when the Reeb vector field $\xi $ belongs to the
distribution ${\frak D}$, we refer to the following

\begin{prop}
\label{prop-0} {\em (Proposition 3.1, \cite{Lee-Suh-2010})} Let $M$ be a
connected orientable Hopf hypersurface in $G_{2}({\Bbb C}^{m+2})$. If the
Reeb vector $\xi $ belongs to the distribution ${\frak D}$, then the
distribution ${\frak D}$ is invariant under the shape operator $A$ of M,
that is $g(A{\frak D},{\frak D}^{\bot })=0$.
\end{prop}

\section{Proof of Theorem~\protect\ref{th-main}\label{sect-main-th}}

In the previous section, because of Lemma \ref{lem-2}, Proposition \ref%
{prop-0} and Theorem \ref{th-A}, we lead to the conclusion that real
hypersurfaces in $G_{2}({\Bbb C}^{m+2})$, under some additional assumptions,
whose normal Jacobi operator is semi-parallel are locally congruent to real
hypersurfaces of type {\bf (A)} or {\bf (B)}. Now, we check if the normal
Jacobi operator of such real hypersurfaces satisfies the semi-parallelism
condition.

\medskip

First, we recall the following proposition due to Berndt and Suh (\cite%
{Berndt-Suh-99}).

\begin{prop}
\label{prop-A} {\em (Proposition 3, \cite{Berndt-Suh-99})} Let $M$ be a
connected real hypersurface of $G_{2}({\Bbb C}^{m+2})$. Suppose that $A%
{\frak D}\subset {\frak D}$, $A\xi =\alpha \xi $ and $\xi $ is tangent to $%
{\frak D}^{\perp }$. Let $J_{1}\in {\frak J}$ be the almost Hermitian
structure such that $JN=J_{1}N$. Then $M$ has three (if $r=\frac{\pi }{2%
\sqrt{8}}$) or four (otherwise) distinct constant principal curvatures
\[
\alpha =\sqrt{8}\cot (\sqrt{8}r),\;\;\beta =\sqrt{2}\cot (\sqrt{2}%
r),\;\;\lambda =-\sqrt{2}\tan (\sqrt{2}r),\;\;\mu =0,
\]%
with some $r\in (0,\frac{\pi }{\sqrt{8}})$. The corresponding multiplicities
are
\[
m(\alpha )=1,\;\;m(\beta )=2,\;\;m(\lambda )=2m-2=m(\mu ),
\]%
and the corresponding eigenspaces are
\begin{eqnarray}
&&T_{\alpha }={\Bbb R}\xi ={\Bbb R}\xi _{1}={\Bbb R}JN={\rm Span}\{\xi \}=%
{\rm Span}\{\xi _{1}\},  \nonumber \\
&&T_{\beta }={\Bbb C}^{\perp }\xi ={\Bbb C}^{\perp }N={\Bbb R}\xi _{2}\oplus
{\Bbb R}\xi _{3}={\rm Span}\{\xi _{2},\xi _{3}\},  \nonumber \\
&&T_{\lambda }=\{X/X\perp {\Bbb H}\xi ,\;\;JX=J_{1}X\},  \nonumber \\
&&T_{\mu }=\{X/X\perp {\Bbb H}\xi ,\;\;JX=-J_{1}X\},\   \nonumber
\end{eqnarray}%
where ${\Bbb R}\xi $, ${\Bbb C}\xi $ and ${\Bbb H}\xi $ respectively denotes
real, complex, quaternionic span of the structure vector field $\xi $ and $%
{\Bbb C}^{\perp }\xi $ denotes the orthogonal complement of the ${\Bbb C}\xi$
in ${\Bbb H}\xi $.
\end{prop}

In this case we have $\xi =\xi _{1}$. From (\ref{eq-normal-Jacobi-operator})
we obtain
\begin{equation}
\overline{R}_{N}(\xi )=8\xi \qquad {\rm and}\qquad \overline{R}_{N}(\xi
_{2})=2\xi _{2}.  \label{eq-B1}
\end{equation}
Since the normal Jacobi operator is semi-parallel, from (\ref%
{eq-semi-parallel}) and the second relation of (\ref{eq-B1}) we obtain:
\begin{equation}
2R(\xi_{2},\xi)\xi_{2} = \overline{R}_{N}(R(\xi_{2},\xi)\xi_{2}),
\label{eq-B2}
\end{equation}
Relation (\ref{eq-Gauss-real-hyp}) for $X=\xi _{2}$, $Y=\xi $ and $Z=\xi
_{2} $ taking into account the fact that $A\xi =\alpha \xi $ and $A\xi
_{2}=\beta \xi _{2}$ implies
\begin{equation}
R(\xi _{2},\xi )\xi _{2}=-(2+\alpha \beta )\xi.  \label{eq-B3}
\end{equation}
Substitution of relation (\ref{eq-B3}) in (\ref{eq-B2}) leads to
\[
(2+\alpha \beta )\xi =0.\
\]
The last relation taking into account that $\alpha=\sqrt{8}\cot(\sqrt{8}r)$
and $\beta=\sqrt{2}\cot(\sqrt{2}r)$ implies
\[
\cot ^{2}(\sqrt{2}r) = 0,
\]%
which is a contradiction. So real hypersurfaces of type {\bf (A)} do not
have semi-parallel normal Jacobi operator.

\medskip

Next we check that whether real hypersurfaces of type {\bf (B)} are equipped
with semi-parallel normal Jacobi operator. We recall the following
proposition due to Berndt and Suh (\cite{Berndt-Suh-99}).

\begin{prop}
\label{prop-B} {\em (Proposition 2, \cite{Berndt-Suh-99})} Let $M$ be a
connected real hypersurface of $G_{2}({\Bbb C}^{m+2})$. Suppose that $A%
{\frak D}\subset {\frak D}$, $A\xi =\alpha \xi $ and $\xi $ is tangent to $%
{\frak D}$. Then the quaternionic dimension $m$ of $G_{2}({\Bbb C}^{m+2})$
is even, say $m=2n$, and $M$ has five distinct constant principal curvatures
\[
\alpha =-2\tan (2r),\;\;\beta =2\cot (2r),\;\;\gamma =0,\;\;\lambda =\cot
(r),\;\;\mu =-\tan (r),
\]%
with some $r\in (0,\pi /4)$. The corresponding multiplicities are
\[
m(\alpha )=1,\;\;m(\beta )=3=m(\gamma ),\;\;m(\lambda )=4n-4=m(\mu ),
\]%
and the corresponding eigenspaces are
\begin{eqnarray}
&&T_{\alpha }={\Bbb R}\xi ={\rm Span}\{\xi \},  \nonumber \\
&&T_{\beta }={\frak J}J\xi ={\rm Span}\{\xi _{1},\xi _{2},\xi _{3}\},
\nonumber \\
&&T_{\gamma }={\frak J}\xi ={\rm Span}\{\varphi _{1}\xi ,\varphi _{2}\xi
,\varphi _{3}\xi \},  \nonumber \\
&&T_{\lambda },\;\;T_{\mu },\   \nonumber
\end{eqnarray}%
where
\[
T_{\lambda }\oplus T_{\mu }=({\Bbb H}{\Bbb C}\xi )^{\perp },\;\;{\frak J}%
T_{\lambda }=T_{\lambda },\;\;{\frak J}T_{\mu }=T_{\mu },\;\;JT_{\lambda
}=T_{\mu }.
\]
\end{prop}

From (\ref{eq-normal-Jacobi-operator}) we obtain
\begin{equation}
\overline{R}_{N}(W)=W,\;\;\;\;\overline{R}_{N}(\xi )=4\xi \;\;{\rm and}\;\;%
\overline{R}_{N}(\xi _{\nu })=4\xi _{\nu },\;\;\nu =1,2,3,  \label{eq-B4}
\end{equation}%
where $W\in T_{\lambda }$. Due to the semi-parallelism of the normal Jaocbi
operator, from (\ref{eq-semi-parallel}) and the first relation of (\ref%
{eq-B4}) we get:
\begin{equation}
R(W,\xi )W=\overline{R}_{N}(R(W,\xi )W),  \label{eq-B5}
\end{equation}%
The Gauss equation (\ref{eq-Gauss-real-hyp}) for $X=W$, $Y=\xi $ and $Z=W$,
because of $A\xi =\alpha \xi $ and $AW=\lambda W$ implies
\begin{equation}
R(W,\xi )W=-\left( 1+\alpha \lambda \right) \xi +\sum_{\nu =1}^{3}g(\varphi
_{\nu }\varphi W,W)\xi _{\nu }.  \label{eq-B6}
\end{equation}%
Substituting (\ref{eq-B6}) in (\ref{eq-B5}) and taking into account relation
(\ref{eq-B4}), we lead to the following

\[
\lbrack 1+\alpha \lambda ]\xi -\sum_{\nu =1}^{3}g(\varphi _{\nu }\varphi
W,W)\xi _{\nu }=0.\
\]%
The inner product of the last relation with $\xi $ and substitution of $%
\alpha =-2\tan (2r)$ and $\lambda =\cot (r)$ yield
\[
1-2\tan (2r)\cot (r)=0,
\]%
from which we obtain
\[
3 + \tan ^{2}(r)=0,
\]%
which is a contradiction. So real hypersurfaces of type {\bf (B)} do not
admit semi-parallel normal Jacobi operator and this completes the proof. $%
\blacksquare $

\bigskip

\noindent {\bf Acknowledgements.} The first author would like to express her
gratitude to Professor Ph. J. Xenos. Second author is thankful to Professor
Oldrich Kowalski for academic hospitality provided by him at Charles
University during June 9-24, 2012.

\bigskip

\noindent Konstantina Panagiotidou

\noindent Mathematics Division-School of Technology,

\noindent Aristotle University of Thessaloniki,

\noindent Thessakibuju 54124, Greece,

\noindent Email: {\tt kapanagi@gen.auth.gr}

\medskip

\noindent Mukut Mani Tripathi

\noindent Department of Mathematics

\noindent Faculty of Science

\noindent Banaras Hindu University

\noindent Varanasi 221005, India

\noindent Email: {\tt mmtripathi66@yahoo.com}


\begin{thebibliography}{99}
\bibitem{Alekseevski-68} D. V. Alekseevskii, {\em Compact quaternion spaces}%
, Func. Anal. Prilo\v{z}en {\bf 2} (1968), no. 2, 11-20.

\bibitem{Berndt-91} J. Berndt, {\em Real hypersurfaces in quaternionic space
forms}, J. Reine Angew. Math. {\bf 419} (1991), 9-26.

\bibitem{Berndt-97} J. Berndt, {\em Riemannian geometry of complex two-plane
Grassmannians}, Rend. Sem. Mat. Univ. Politec. Torino {\bf 55} (1997), no.
1, 19-83.

\bibitem{Berndt-Suh-99} J. Berndt and Y. J. Suh, {\em Real hypersurfaces in
complex two-plane Grassmannians}, Monatsh. f\"{u}r Math. {\bf 127} (1999),
no. 1, 1-14.

\bibitem{Berndt-Suh-02} J. Berndt and Y. J. Suh, {\em Real hypersurfaces
with isometric Reeb flow in complex two-plane Grassmannians}, Monatsh. f\"{u}%
r Math. {\bf 137} (2002), no. 2, 87-98.

\bibitem{CK} J. T. Cho and M. Kimura, {\em Curvature of Hopf hypersurfaces
in a complex space form}, Results Math. {\bf 61} (2012), no. 1-2, 127--135.

\bibitem{Jeong-Suh-2008} I. Jeong and Y. J. Suh, {\em Real hypersurfaces in
complex two-plane Grassamnnians with Lie $\xi $-parallel normal Jacobi
operator}, J. Korean Math. Soc. {\bf 45} (2008), no. 4, 1113-1133.

\bibitem{Jeong-Kim-Suh-2010} I. Jeong, H. J. Kim and Y. J. Suh, {\em Real
hypersurfaces in complex two-plane Grassamnnians with parallel normal Jacobi
operator}, Publ. Math. Debrecen {\bf 76} (2010), no.1-2, 203-218.

\bibitem{Jeong-Lee-Suh-2011} I. Jeong, H. Lee and Y. J. Suh, {\em Hopf
hypersurfaces in complex two-plane Grassamnnians with Lie parallel normal
Jacobi operator}, Bull. Korean Math. Soc. {\bf 48} (2011), no. 2, 427-444.

\bibitem{Jeong-Suh-2011} I. Jeong and Y. J. Suh, {\em Real hypersurfaces in
complex two-plane Grassamnnians with ${\frak F}$-parallel normal Jacobi
operator}, Kyungpook Math. J. {\bf 51} (2011) no. 4, 395-410.

\bibitem{Jeong-Suh-Tripathi-2012} I. Jeong, Y. J. Suh and M. M. Tripathi,
{\em Real hypersurfaces of type $A$ in complex two-plane Grassamnnians
related to the normal Jacobi operator}, Bull. Korean Math. Soc. {\bf 49}
(2012).

\bibitem{Lee-Suh-2010} H. Lee and Y. J. Suh, {\em Real hypersurfaces of type
B in complex two-plane Grassamnnians related to the Reeb vector}, Bull.
Korean Math. Soc. {\bf 47} (2010), no. 3, 551-561.

\bibitem{Machado-Perez-Suh-2011} C. J. G. Machado, J. D. Perez, I. Jeong and
Y. J. Suh, {\em Real hypersurfaces in complex two-plane Grassamnnians whose
normal Jacobi operator is of Codazzi type}, Cent. Eur. J. Math. {\bf 9}
(2011), no. 3, 578-582.

\bibitem{NR1} R. Niebergall and P. J. Ryan, {\em Semi-parallel and
semi-symmetric real hypersurfaces in complex space forms}, Kyungpook Math.
J. {\bf 38} (1998), 227-234.

\bibitem{Perez-Jeong-Suh-2007} J. D. Perez, I. Jeong and Y. J. Suh, {\em %
Real hypersurfaces in complex two-plane Grassamnnians with commuting normal
Jacobi operator}, Acta Math. Hungar. {\bf 117} (2007), no. 3, 201-217

\bibitem{PS1} J. D. Perez and F. G. Santos, {\em Real hypersurfaces in
complex projective space whose structure Jacobi operator is cyclic-Ryan
parallel}, Kyungpook Math. J. {\bf 49} (2009), 211-219.
\end{thebibliography}
\end{document}